\def\draft{n}
\documentclass{amsart}
\usepackage[headings]{fullpage}
\usepackage{amssymb,epic,eepic,epsfig}

\theoremstyle{plain}

\newtheorem{theorem}{Theorem}
\newtheorem{proposition}{Proposition}[section]
\newtheorem{lemma}[proposition]{Lemma}

\theoremstyle{definition}
\newtheorem{definition}[proposition]{Definition}

\newtheorem{question}{Question}

\theoremstyle{remark}

\newtheorem{remark}[proposition]{Remark}

\def\printname#1{
        \if\draft y
                \smash{\makebox[0pt]{\hspace{-0.5in}
                        \raisebox{8pt}{\tt\tiny #1}}}
        \fi
}

\newlength{\standardunitlength}
\catcode`\@=11
\long\def\@makecaption#1#2{%
     \vskip 10pt

\setbox\@tempboxa\hbox{
       \small\sf{\bfcaptionfont #1. }\ignorespaces #2}%
     \ifdim \wd\@tempboxa >\captionwidth {%
         \rightskip=\@captionmargin\leftskip=\@captionmargin
         \unhbox\@tempboxa\par}%
       \else
         \hbox to\hsize{\hfil\box\@tempboxa\hfil}%
     \fi}
\font\bfcaptionfont=cmssbx10 scaled \magstephalf
\newdimen\@captionmargin\@captionmargin=2\parindent
\newdimen\captionwidth\captionwidth=\hsize
\catcode`\@=12

\def\lbl#1{\label{#1}\printname{#1}}

\def\BN{\mathbb N}
\def\BZ{\mathbb Z}

\def\BQ{\mathbb Q}

\def\BC{\mathbb C}

\def\a{\alpha}

\def\La{\Lambda}
\def\l{\lambda}
\def\Ga{\Gamma}

\def\la{\langle}
\def\ra{\rangle}

\def\e{\epsilon}
\def\Ga{\Gamma}

\def\b{\beta}

\def\ft{\mathrm{ft}}
\def\gterm{term}
\def\bterm{balanced term}
\def\ft{\mathfrak{t}}
\def\lcm{\mathrm{lcm}}
\def\ord{\mathrm{ord}}
\def\supp{\mathrm{supp}}

\begin{document}


\title[$G$-functions and multisum versus holonomic sequences]{
$G$-functions and multisum versus holonomic sequences}
\author{Stavros Garoufalidis}
\address{School of Mathematics \\
         Georgia Institute of Technology \\
         Atlanta, GA 30332-0160, USA \\ 
         {\tt http://www.math.gatech} \newline {\tt .edu/$\sim$stavros } }
\email{stavros@math.gatech.edu}

\thanks{The author was supported in part by NSF. \\
\newline
1991 {\em Mathematics Classification.} Primary 57N10. Secondary 57M25.
}

\date{November 11, 2008 }


\begin{abstract}
The purpose of the paper is three-fold: (a) we prove that every sequence
which is a multidimensional sum of a balanced hypergeometric term has
an asymptotic expansion of Gevrey type-1 with rational exponents, (b)
we construct a class of $G$-functions that come from enumerative
combinatorics, and (c) we give a counterexample to a question of Zeilberger
that asks whether holonomic sequences can be written as multisums of
balanced hypergeometric terms. The proofs utilize the notion of a $G$-function,
introduced by Siegel, and its analytic/arithmetic properties shown recently 
by Andr\'e.
\end{abstract}

\maketitle

\tableofcontents

\section{Introduction}
\lbl{sec.intro}

\subsection{Balanced multisum sequences are of Nilsson type}
\lbl{sub.balanced}

The purpose of the paper is three-fold: 

\begin{itemize}
\item[(a)] we prove that every sequence
which is a multidimensional sum of a balanced hypergeometric term has
an asymptotic expansion of Gevrey type-1 with rational exponents, 
\item[(b)]
we construct a class of $G$-functions that come from enumerative
combinatorics, and 
\item[(c)] we give a counterexample to a question of Zeilberger
that asks whether holonomic sequences can be written as multisums of
balanced hypergeometric terms. 
\end{itemize}
The proofs utilize the notion of a $G$-function, introduced by Siegel, 
and its analytic/arithmetic properties shown recently by Andr\'e.
Let us begin by introducing the notion of a (balanced) multisum sequence.

\begin{definition}
\lbl{def.bc}
A {\em (balanced) multisum sequence} $(a_n)$ is a sequence of complex numbers
of the form
\begin{equation}
\lbl{eq.bc}
a_n=\sum_{k \in \supp(\ft_{n,\bullet})} \ft_{n,k}
\end{equation}
where $\ft$ is a (balanced) term and the sum is over a finite set that
depends on $\ft$.
\end{definition}

\begin{definition}
\lbl{def.hyperg}
A {\em term} $\ft_{n,k}$ in variables $(n,k)$ where $k=(k_1,\dots,k_r)$
is an expression of the form:
\begin{equation}
\lbl{eq.defterm}
\ft_{n,k}=C_0^n \prod_{i=1}^r C_i^{k_i} \prod_{j=1}^J A_j(n,k)!^{\e_j}
\end{equation}
where $C_i \in \overline{\BQ}$ for $i=0,\dots,r$, 
$\e_j=\pm 1$ for $j=1,\dots,J$, and $A_j$ are integral 
linear forms in the variables $(n,k)$ such that for every 
$n \in \BN$, the set 
\begin{equation}
\lbl{eq.kset}
\supp(\ft_{n,\bullet}):=
\{ k \in \BZ^r \, | \, A_j(n,k) \geq 0, \,\, j=1,\dots, J \}
\end{equation}
is finite. We will call a \gterm\ {\em balanced} if in addition it
satisfies the {\em balance condition}:
\begin{equation}
\lbl{eq.Ajsum}
\sum_{j=1}^J \e_j A_j=0.
\end{equation}
\end{definition}
For example, the {\em Ap\'ery sequence} (see \cite{vdP})
is a balanced multisum given by:
\begin{equation}
\lbl{eq.apery}
a_n=\sum_{k=0}^n \binom{n}{k}^2 \binom{n+k}{k}^2=
\sum_{k=0}^n \left(\frac{(n+k)!}{k!^2 (n-k)!}\right)^2.
\end{equation}
Multisum sequences appear frequently in enumeration questions; for numerous
examples, see \cite{St2,FS}. A key problem is to study the asymptotics of a 
(balanced) multisum sequence. This is a classical problem that has been
discussed by several authors for decades, see \cite{BT,FS,St2,WP}. Parsing
through the literature, in numerous examples of balanced multisum examples, 
a certain rationality of the 
leading exponents of $n$ was found by accident, with no explanation. 
Understanding this rationality lead to the results of our paper.

To explain this rationality, let us introduce one more
definition.

\begin{definition}
\lbl{def.nilsson}
We say that a sequence $(a_n)$ is of {\em Nilsson type} if it has an 
asymptotic expansion of the form
\begin{equation}
\lbl{eq.anexp}
a_n \sim \sum_{\l,\a,\b} \l^{-n} n^{\a} (\log(n))^{\b} 
f_{\l,\a,\b}\left(\frac{1}{n}\right)
\end{equation}
where the summation is over a finite set of triples $(\l,\a,\b)$,
$\l \in \overline\BQ$ is an algebraic number, $\a \in \BQ$, $\b \in \BN$,
and $f_{\l,\a,\b}(z)$ is Gevrey-1 power series, i.e., the coefficient of 
$z^k$ in $f_{\l,\a,\b}(z)$ is bounded by $C^n$ for some $C>0$.
\end{definition}
Now we can state our first result which does not seem to be covered by the
existing literature on asymptotic expansions of sequences.

\begin{theorem}
\lbl{thm.1}
Every balanced multisum sequence is of Nilsson type.
\end{theorem}

\subsection{$G$-functions}
\lbl{sub.Gfunctions}

The proof of Theorem \ref{thm.1} utilizes the notion of a $G$-function,
introduced by Siegel in \cite{Si} with motivation being
arithmetic problems in elliptic integrals, and transcendence
problems in number theory. For further information about $G$-functions
and their properties, see \cite{An,Bo,DGS,Si,To}.

\begin{definition}
\lbl{def.Gfunction}
We say that series $G(z)=\sum_{n=0}^\infty a_n z^n$ is a {\em $G$-function}
if 
\begin{itemize}
\item[(a)]
the coefficients $a_n$ are algebraic numbers,
\item[(b)]
there exists a constant $C>0$ so that for every $n \in \BN$
the absolute value of every conjugate of $a_n$ is less than or equal to 
$C^n$, 
\item[(c)]
the common denominator of $a_0,\dots, a_n$ is less than or equal 
to $C^n$,
\item[(d)]
$G(z)$ is holonomic, i.e., it satisfies a linear differential equation
with coefficients polynomials in $z$.
\end{itemize}
\end{definition}

\begin{remark}
\lbl{rem.ag0}
In \cite{An}, Andr\'e calls a series that satisfies (a)-(c), an {\em
arithmetic Gevrey-0} series.  
\end{remark}

Our next theorem is a construction of $G$-functions from balanced terms.
Notice that it is easy to generate examples of balanced terms; see for example
the Ap\`ery sequence above.

\begin{theorem}
\lbl{thm.2}
If $(a_n)$ is a balanced multisum sequence, the generating series
$G(z)=\sum_{n=0}^\infty a_n z^n$ is a $G$-function.
\end{theorem}
$G$-functions can be constructed by arithmetic, or by geometry. Theorem
\ref{thm.2} offers a new construction.
For further discussion, see Section \ref{sec.discuss}.

The proof of Theorem \ref{thm.2} also gives the following result, 
which may be of interest to enumerative combinatorial problems that often
lead to sequences of integers.

\begin{theorem}
\lbl{thm.22}
The generating series of an integer-valued, exponentially bounded 
holonomic sequence is a $G$-function.
\end{theorem}

\subsection{Holonomic sequences are not balanced multisums}
\lbl{sub.holo}

We now come to the third part of the paper, which compares multisum and
holonomic sequences. 

\begin{definition}
\lbl{def.holo}
A sequence $(a_n)$ is
{\em holonomic} (i.e., $D$-{\em finite} in the sense of \cite{St1}) 
if it satisfies a linear recursion relation with polynomial coefficients.
In other words, there exist $d \in \BN$ and $P_j(n) \in \overline{\BQ}[n]$ 
(where $\overline{\BQ}$ denotes the set of algebraic numbers) for
$j=1,\dots,d$, so that for every $n \in \BN$ we have:
\begin{equation}
\lbl{eq.Prec}
P_d(n) a_{n+d} + \dots P_0(n) a_n =0.
\end{equation}
\end{definition}

The following is a fundamental theorem of Wilf-Zeilberger.

\begin{theorem}
\lbl{thm.z}\cite{WZ,Ze}
Every multisum sequence is holonomic.
\end{theorem}

The above theorem has a constructive proof with several computer 
implementations, see \cite{PR, PWZ} and \cite{WZ}. The converse 
was widely accepted as a reasonable conjecture, communicated
to the author by Zeilberger. Our goal is to give a counterexample, and give 
an obstruction for the converse to hold. 

\begin{theorem}
\lbl{thm.3}
Consider the holonomic sequence $(a_n)$ defined by 
\begin{equation}
\lbl{eq.counter}
(2n+1)a_{n+2}-(7n+11)a_{n+1}+(2n+1)a_n=0.
\end{equation}
with initial conditions $a_0=0$, $a_1=1$. Then, $(a_n)$ is not a balanced
multisum.
\end{theorem}
To understand why the converse to Theorem \ref{thm.z} fails, and why
the example given by \eqref{eq.counter} is not pathological 
(but rather typical), let us look at the asymptotic expansion of an
exponentially bounded 
holonomic sequence. It follows from Birkhoff-Trjitzinsky and Turrittin
that a holonomic sequence is almost of Nilsson type, i.e., it satisfies
\eqref{eq.anexp} where the exponents $\a$ are algebraic, but not necessarily
rational numbers. This also typically happens in the analysis of linear
ODE, where the above exponents are known as {\em Frobenius exponents}.
For more details, see Section \ref{sub.thm3}. On the other hand, balanced 
multisum sequences have rational exponents according to Theorem \ref{thm.1}.
This proves and explains Theorem \ref{thm.3}.

\subsection{Acknowledgement}
The author wishes to thank Y. Andr\'e and D. Zeilberger 
for enlightening conversations, guidance and encouragement.

\section{Proofs}
\lbl{sec.proofs}

\subsection{Proof of Theorem \ref{thm.2}}
\lbl{sub.thm2}

Let us begin with the following alternative presentation of a \bterm .

\begin{lemma}
\lbl{lem.bterm}
Every \bterm\ $\ft$ can be written in the form:
\begin{equation}
\lbl{eq.bterm}
\ft_{n,k}=C_0^n \prod_{i=1}^r C_i^{k_i} \prod_{j=1}^J 
\binom{B_j(n,k)}{D_j(n,k)}^{\e_j}
\end{equation}
where $C_i \in \overline{\BQ}$ for $i=0,\dots,r$, 
$\e_j=\pm 1$ for $j=1,\dots,J$, and $B_j, D_j$ are integral 
linear forms in the variables $(n,k)$. 
\end{lemma}

\begin{proof}
Consider a \bterm\ $\ft$ given by \eqref{eq.defterm}, where the linear
forms $A_j$ satisfy the balance condition \eqref{eq.Ajsum}. Let 
$J^{\pm}=\{j \in J \, | \e_j=\pm 1\}$ and consider 
the linear form $A(n,k)$ defined by:
$$
A(n,k)=\sum_{j \in J^+} A_j(n,k)=\sum_{j \in J^-} A_j(n,k),
$$
where the second equality follows from the balance condition. Then,
multiply and divide the \bterm\ by $A(n,k)!$, and rearrange the factors
into a ratio of multibinomial coefficients as follows:

\begin{eqnarray*}
\prod_{j=1}^J (A_j(n,k)!)^{\e_j} &=& 
\frac{\prod_{j \in J^+} A_j(n,k)!}{\prod_{j \in J^-} A_j(n,k)!} \\
&=& 
\frac{\prod_{j \in J^+} A_j(n,k)!}{\prod_{j \in J^-} A_j(n,k)!} 
\frac{A(n,k)!}{A(n,k)!} \\
&=& \frac{\binom{A(n,k)}{A_j \, | j \in J^-}}{
\binom{A(n,k)}{A_j \, | j \in J^+}}
\end{eqnarray*}
Now, write the multibinomial coefficient as a product of binomial coefficients.
The result follows.
\end{proof}

The next lemma from number theory is well-known (see \cite[p.198]{vdP}
and also \cite{Si}) and follows easily from Chebytchev's theorem. 
Below, lcm denotes the {\em least common multiple}.


\begin{lemma}
\lbl{lem.lcm}\cite{vdP,Si}
There exists $C>0$ so that
\begin{equation}
\lbl{eq.lcm}
\lcm\left(\binom{n}{0},\dots, \binom{n}{n}\right) < C^n
\end{equation}
for all $n \in \BN$.
\end{lemma}

\begin{proof}
Let $\ord_p m $ denote the maximal power of a prime number $p$ that divides a 
natural number $m$. Then, for every natural number $n$ and integer $a$ and 
$b$ with $0 \leq b \leq a \leq n$ and 
every prime number $p$ we have:
$$
\ord_p \binom{a}{b} \leq \left[\frac{\log a}{\log p} \right] -\ord_p b \leq
\ord_p \lcm(1,\dots,a) -\ord_p b \leq \ord_p \lcm(1,\dots,n).
$$
Thus,
$$
\lcm\left(\binom{n}{0},\dots, \binom{n}{n}\right)   \leq  \lcm(1,\dots,n)
$$
On the other hand, it is known that
$$
\log \lcm(1,\dots,n) =O(n).
$$
For a detailed discussion, see \cite[p.198]{vdP}.
\end{proof}

We are now ready to give the proof of Theorem \ref{thm.2}.

\begin{proof}(of Theorem \ref{thm.2})
Fix a \bterm\ $\ft_{n,k}$ as in \eqref{eq.bterm}, and the corresponding 
sequence $(a_n)$ of \eqref{eq.bc}. 
We will show that conditions (a),(b),(c) of Definition 
\ref{def.Gfunction} are satisfied. Condition (a) is obvious.

Using $\binom{a}{b} \leq 2^a$, and the fact that the set \eqref{eq.kset} is
a subset of $[-K n, K n]^r \cap \BZ^r$ for some $K>0$, Equation 
\eqref{eq.bterm} implies that there exists a constant $C >0$ so that
$$
|\ft_{n,k}| < C^n
$$
for all $(n,k)$ and for all complex conjugates of $\ft_{n,k}$. Summing up with
respect to $k$ in Equation \eqref{eq.bc}, and using the fact that the 
summation set has polynomial size in $n$, it follows (after possibly enlarging
$C$) that 
$$
|a_{\ft,n}| < C^n
$$
for all $n >0$ and for all complex conjugates of $a_n$. This proves
condition (b) of Definition \ref{def.Gfunction}. 

Condition (c) follows from
Equation \eqref{eq.bterm}, Lemma \ref{lem.lcm} and the fact that the summation
set \eqref{eq.kset} is bounded polynomially by $n$.

Condition (d) follows from Wilf-Zeilberger's Theorem \ref{thm.z}.
\end{proof}

\subsection{The local monodromy of a $G$-function}
\lbl{sub.Gmono}

In this section we will make little distinction between a convergent
power series, its analytic
continuation, and the corresponding function. Recall that a power series
is {\em holonomic} if it satisfies a linear differential equation 
$P G(z)=0$ where $P \in \overline{\BQ}\la z, d/dz \ra$ is a linear differential
operator with coefficients in $\overline\BQ[z]$.
By the theory of differential
equations (see for example \cite{In,O}), 
a holonomic function has analytic continuation as a multivalued
analytic function in $\BC\setminus\La$, where $\La$ is a finite set of
algebraic numbers. The following Theorem follows from a combination
of results of Katz, Andr\'e and Chudnovsky; see \cite{Ka,An,CC} and also
\cite{C-L} for a detailed exposition.

\begin{theorem}
\lbl{thm.andre}\cite{An,Ka,CC}
The local monodromy $T$ of a $G$-function around a singularity is
quasi-unipotent. In other words, 
\begin{equation}
\lbl{eq.Tmonodromy}
(T^r-1)^s=0
\end{equation} 
for some nonzero natural numbers $r$ and $s$. 
\end{theorem}
It follows that the local expansion of $G(z)$ around a singularity 
$\l \in \La$ is a finite sum of series of the form:
\begin{equation}
\lbl{eq.nilsson}
\sum_{\a,\b} c_{\a,\b} (z-\l)^{\a} (\log(z-\l))^{\b} h_{\a,\b}(z-\l)
\end{equation}
where $\a \in \BQ$, $\b \in \BN$, $c_{\a,\b} \in \BC$
and $h_{\a,\b}(w)$ are convergent germs
at $w=0$. In fact, Andr\'e shows that $h_{\a,\b}(w)$ are $G$-functions 
themselves. Said differently, the $G$-functions that come from arithmetic 

\begin{itemize}
\item[(a)]
are {\em regular holonomic} (i.e., the power series $h_{\a,\b}(w)$ above are
convergent at $w=0$), and 
\item[(b)]
have 
{\em rational exponents} (denoted by $\{\a\}$ above).
\end{itemize}
On the other hand, the generating series of a generic exponentially bounded
holonomic sequence $(a_n)$ will not be in general regular holonomic,
nor will they have rational exponents. This explains Theorem \ref{thm.3}.

\begin{remark}
Power series of the form \eqref{eq.nilsson} are 
known in the literature as {\em Nilsson series}; see \cite{Ni}.
\end{remark}

\subsection{The Taylor series of a $G$-function and Theorem \ref{thm.1}}
\lbl{sub.propar}

The following lemma is a well-known application of Cauchy's theorem; see
for example \cite[Thm.A]{Ju}.

\begin{lemma}
\lbl{lem.jungen}\cite[Thm.A]{Ju}
If $\a \in \BC\setminus\BN$, $\b \in \BN$, and 
$$
(1-z)^{\a} (\log(1-z))^{\b}=\sum_{n=0} a_n z^n
$$
then
$$
a_n=\frac{n^{-\a-1}}{\Ga(-\a)}\left((\log(n))^{\b} \phi_0(n) + \dots
(\log(n))^{0} \phi_{\b}(n)\right)
$$
where $\phi_j(z)$ for $j=0,\dots,\b$
are Gevrey-$1$ series with rational coefficients.
\end{lemma}

Recall that a series $\sum_{n=0}^\infty a_n z^n$ is {\em Gevrey}-1 (resp. 
{\em arithmetic Gevrey}-1) 
if $\sum_{n=0}^\infty (a_n/n!) z^n$ is convergent at $z=0$ 
(resp. a $G$-function). Lemma \ref{lem.jungen} and a deformation of the
contour argument implies the following. See also \cite[Sec.7]{CG}.

\begin{proposition}
\lbl{prop.expar}
If $G(z)=\sum_{n=0}^\infty a_n z^n$ is a $G$-function, then
$(a_n)$ is of Nilsson type, where $\l$ in Equation \eqref{eq.anexp}
are the singularities of $G(z)$, and $\a$, $\b$, $f_{\l,\a,\b}$ 
are determined by the local monodromy of $G(z)$ at $z=\l$.
\end{proposition}

\begin{remark}
\lbl{rem.expar}
In fact, Jungen's proof combined with Andr\'e's theorem that the series
$h_{\a,\b}(w)$ are $G$-functions, implies that the series $f_{\l,\a,\b}(z)$
of Equation \eqref{eq.anexp} are arithmetic Gevrey-1. 
\end{remark}

\subsection{Proof of Theorem \ref{thm.3}}
\lbl{sub.thm3}

The next lemma is well-known. 

\begin{lemma}
\lbl{lem.1}
If $(a_n)$ is holonomic, the generating series $G(z)=\sum_{n=0}^\infty a_n z^n$
is holonomic.
\end{lemma}
Birkhoff-Trjitzinsky, followed by Turrittin 
(see \cite{BT,Tu,Pr} and \cite[Eqn.1.3]{BC}) prove the following result
concerning the asymptotic expansion of a holonomic sequence.

\begin{proposition}
\lbl{prop.exphol}
If $G(z)=\sum_{n=0}^\infty a_n z^n$ is a holonomic function, then
\begin{equation}
\lbl{eq.anexphol}
a_n \sim \sum_{\l,\a,\b,s} n!^s \l^{-n} n^{-\a-1} (\log(n))^{\b}
f_{\l,\a,\b,s}\left(\frac{1}{n}\right)
\end{equation}
where $\l$ lies in a subset of the finite set of singularities of $G(z)$, 
$s$ lies in a finite set of nonpositive rational numbers, and 
$\a, \b$ are the exponents in the local expansion of $G(z)$
around $\l$, and $f_{\l,\a,\b,s}(z)$ are Gevrey-1.
\end{proposition}

\subsection{The exponents of the sequence of Theorem \ref{thm.1}}
\lbl{sub.thm1}


It remains to compute the exponents of the holonomic function $G(z)$
associated to the sequence of Equation \eqref{eq.counter}. 
One way to solve this problem is to 
convert the holonomic equation \eqref{eq.counter} into a 
differential equation for the generating series
and compute the exponents of the differential equation using {\em Frobenius's
method}; see \cite{O,In}. In addition, one needs to show that the corresponding
constants $c_{\a,\b}$ in \eqref{eq.nilsson} are nonvanishing.
An alternative way is to relate the exponents of the generating series $G(z)$
of a sequence $(a_n)$ to the asymptotic expansion of the sequence itself.

Consider the sequence $(a_n)$ given by \eqref{eq.counter} and its generating
series $G(z)$. Converting the recursion relation for $(a_n)$ into a 
differential equation for $G(z)$ we obtain that $G(z)$ 
satisfies the inhomogeneous differential equation:

\begin{equation}
\lbl{eq.GODE}
z(z^2 -7z+2)G'(z) +(z^2-4z-3)G(z)+z=0, \qquad G(0)=0
\end{equation}
If we wish, we can divide by $z$ and differentiate once to get a linear
second order differential equation for $G(z)$. The singularities $\La$ 
of $G(z)$ is a subset of the roots of $z(z^2 -7z+2)$. I.e., we have:

\begin{equation}
\lbl{eq.lambda}
\La \subset \{0, \frac{1}{4}(7 \pm \sqrt{33}) \}.
\end{equation}
Frobenius's method gives that the exponent at 
$\l_{\pm}=\frac{1}{4}(7 \pm \sqrt{33})$ is given by

\begin{equation}
\lbl{eq.expa}
\a_{\pm}=-1\pm \frac{5}{2} \sqrt{\frac{3}{11}}
\end{equation}
which is non-rational. It is easy to compute that $\b=0$. It remains to
argue that the so-called {\em Stokes constant} $c_{\a_{\pm},\b} \neq 0$. 
One can do an explicit numerical computation in the spirit of 
\cite[Sec.4]{FT}, using {\em Pad\'e approximants} and working in the so-called
{\em Borel plane}.

Alternatively, we may argue as follows. If $G(z)$ is analytic at 
$\frac{1}{4}(7 + \sqrt{33})$, then by Galois invariance and Equation 
\eqref{eq.lambda}, it follows that $G(z)$ is entire. Lemma \ref{lem.entireG}
below implies that $G(z)$ is a polynomial. It follows that $a_n=0$ for 
sufficiently large $n$. The recursion relation \eqref{eq.counter} implies
that $a_n=0$ for all $n \in \BN$, an obvious contradiction.
\qed

The next lemma was communicated to us by Y. Andr\'e, and is a useful way of
deducing the existence of singularities of $G$-functions. For a detailed
discussion, see also \cite{C-L}.

\begin{lemma}
\lbl{lem.entireG}
Every entire $G$-function is a polynomial.
\end{lemma}

\begin{proof}
According to Chudnovsky and Katz, a $G$-function is a solution of a {\em
Fuchsian} differential equation, i.e., regular singular in 
$\BC \cup\{\infty\}$. An entire $G$-function does not have
any monodromy at finite distance, hence it does not have any monodromy
at infinity as well. According to a classical result of Schlesinger, any
solution of a Fuchsian differential equation which is invariant under
the global monodromy group is a rational function. If, moreover, it is
entire, then it is a polynomial function.
\end{proof}

\section{Further discussion}
\lbl{sec.discuss}

Theorem \ref{thm.2} may be viewed as a way of constructing holonomic 
$G$-functions from enumerative combinatorics.
There are two well-known sources of $G$-functions: from
{\em arithmetic} (see Theorem \ref{thm.andre} and also \cite{An,Bo,DGS,Ka}), 
and from {\em geometry}, related to the regularity of the 
{\em Gauss-Manin connection}. For the latter, see for example
example, \cite{Br,De,Ka}. In all cases (combinatorics, geometry and 
arithmetic), the constructed $G$-functions are regular holonomic with
rational exponents.

The $G$-functions obtained geometry and arithmetic are closely 
related. The main conjecture is that all $G$-functions come from
geometry. For a discussion of this topic, and for a precise formulation
of the Bombieri-Dwork Conjecture, see the survey papers of
\cite{Bo,Ka} and also \cite[p.8]{To}.
Our question is motivated by Theorem \ref{thm.2} and Bombieri-Dwork
Conjecture of \cite[p.8]{To}.

\begin{question}
\lbl{que.1}
If $(a_n)$ is an integer valued, exponentially bounded holonomic sequence,
does it follow that it is a multisum sequence?
\end{question}

Our next question compares the $G$-functions of Theorem \ref{thm.2} with
those that come from geometry.

\begin{question}
\lbl{que.2}
Does every $G$-function of Theorem \ref{thm.2} come from geometry?
\end{question}
In \cite{Ga1} this was shown to be true when the balanced term $\ft_{n,k}$ is
{\em special}, i.e., it is a product of binomials of linear forms of $(n,k)$
(in other words, $\e_j=+1$ for all $j=1,\dots,J$ in Equation \eqref{eq.bterm}.

Finally, let us point out that the proof of Theorem \ref{thm.1} and
Theorem \ref{thm.2} is not constructive. In particular, it would be nice
to be able to compute the singularities of the generating series of
a balanced multisum sequence directly from the \bterm $\ft$. With this
in mind, the author 
developed an efficient ansatz for the asymptotics of balanced multisum
sequences; see \cite{Ga1}. When $r=1$ in Equation \eqref{eq.defterm} (i.e.,
for single-sums), the ansatz can be proven using the Euler-MacLaurin formula
and various ideas of resurgence; see \cite{Ga2}.

\ifx\undefined\bysame
        \newcommand{\bysame}{\leavevmode\hbox
to3em{\hrulefill}\,}
\fi


\begin{thebibliography}{[EMSS]}

\bibitem[An]{An} Y. Andr\'e, 
        {\em S\'eries Gevrey de type arithm\'etique. I. Th\'eor\`emes 
        de puret\'e et de dualit\'e}, 
        Ann. of Math. (2)  {\bf 151}  (2000) 705--740.

\bibitem[BC]{BC} M.A. Barkatou and G. Chen, 
        {\em Some formal invariants of linear difference systems and their 
        computations}, 
        J. Reine Angew. Math. {\bf 533} (2001) 1--23. 

\bibitem[BT]{BT} G. Birkhoff and W. Trjitzinsky,
        {\em Analytic theory of singular difference equations},
        Acta Math. {\bf 60} (1932) 1--89.
 
\bibitem[Bo]{Bo} E. Bombieri, 
        {\em On $G$-functions},
        in Recent progress in analytic number theory, Academic Press
        Vol. 2 (1981) 1--67.

\bibitem[Br]{Br} E. Brieskorn,
        {\em Die Monodromie der isolierten Singularit\"aten von 
        Hyperfl\"achen},  
        Manuscripta Math. {\bf 2}  (1970) 103--161.

\bibitem[C-L]{C-L} A. Chambert-Loir, 
        {\em Th\'eor\`emes d'alg\'ebricit\'e en g\'eom\'etrie diophantienne 
        (d'apr\`es J.-B. Bost, Y. Andr\'e, D. \& G. Chudnovsky)},
        S\'eminaire Bourbaki, Vol. 2000/2001. 
        Ast\'erisque No. {\bf 282} (2002) Exp. No. 886  175--209. 

\bibitem[CC]{CC} D.V. Chudnovsky and G.V. Chudnovsky, 
        {\em Applications of Pad\'e approximations to the Grothendieck 
        conjecture on linear differential equations},
        in  Number theory  Lecture Notes in Math. {\bf 1135} 
        Springer-Verlag (1985) 52--100.

\bibitem[CG]{CG} O. Costin and S. Garoufalidis,
        {\em Resurgence of the Kontsevich-Zagier power series},
        Annales de l' Institut Fourier, {\em in press}.

\bibitem[DGS]{DGS} B. Dwork, G. Gerotto and F.J. Sullivan, 
        {\em An introduction to $G$-functions},
        Annals of Mathematics Studies, Princeton University Press 
        {\bf 133} 1994. 

\bibitem[De]{De} P. Deligne,
        {\em \'Equations diff\'erentielles \`a points singuliers r\'eguliers},
        Lecture Notes in Mathematics, {\bf 163} Springer-Verlag 1970.

\bibitem[FT]{FT} F. Fauvet and J. Thomann, 
        {\em Formal and numerical computations with resurgent functions},
        Numer. Algorithms {\bf 40} (2005) 323--353.
 
\bibitem[FS]{FS} P. Flajolet and R. Sedgewick,
        {Analytic combinatorics}, Cambridge University Press 2008.

\bibitem[Ga1]{Ga1} S. Garoufalidis,
        {\em An ansatz for the asymptotics of hypergeometric 
        multisums},
        Advances in Applied Mathematics, {\bf 41} (2008) 423--451.

\bibitem[Ga2]{Ga2} \bysame,
        {\em Resurgence of 1-dimensional hypergeometric multisums}, to appear.

\bibitem[In]{In} E.L. Ince, 
        {\em Ordinary Differential Equations}, 
        Dover Publications, 1944.
 
\bibitem[Ju]{Ju} R. Jungen,
        {\em Sur les s\'eries de Taylor n'ayant que des singularit\'es 
        alg\'ebrico-logarithmiques sur leur cercle de convergence},
        Comment. Math. Helv. {\bf 3} (1931) 266--306. 

\bibitem[Ka]{Ka} N.M. Katz,
        {\em Nilpotent connections and the monodromy theorem: Applications 
        of a result of Turrittin},
        IHES Publ. Math. {\bf 39} (1970) 175--232. 

\bibitem[Ni]{Ni} N. Nilsson, 
        {\em Some growth and ramification properties of certain integrals 
        on algebraic manifolds},
        Ark. Mat. {\bf 5} (1965) 463--476.

\bibitem[O]{O} F. Olver, 
        {\em Asymptotics and special functions}, Reprint. 
        AKP Classics. A K Peters, Ltd., Wellesley, MA, 1997.
 
\bibitem[PR]{PR} P. Paule and A. Riese, Mathematica software:
        \verb+http://www.risc.uni-linz.ac.at/research/combinat/risc/software/+

\bibitem[PWZ]{PWZ} M. Petkov\v sek, H.S. Wilf and D.Zeilberger,
        {\em $A=B$},
        A.K. Peters, Ltd., Wellesley, MA 1996.

\bibitem[vdP]{vdP} A. van der Poorten, 
        {\em A proof that Euler missed$\ldots $Ap\`ery's proof of the 
        irrationality of $\zeta (3)$},
        Math. Intelligencer {\bf 1} (1978/79) 195--203. 

\bibitem[Pr]{Pr} C. Praagman, 
        {\em The formal classification of linear difference operators},
        Nederl. Akad. Wetensch. Indag. Math. {\bf 45} (1983) 249--261. 

\bibitem[Si]{Si} C.L. Siegel,
        {\em \"Uber einige Anwendungen diophantischer Approximationen},
        Abh. Preuss. Akad. Wiss. {\bf 1} (1929) 1--70. Reprinted in 
        Gesammelte Abhandlungen, vol. 1, no 16 (1966) 209--266. 

\bibitem[St1]{St1} R.P. Stanley,
        {\em Differentiably finite power series},
        European J. Combin. {\bf 1} (1980) 175--188. 

\bibitem[St2]{St2} \bysame,
        {\em Enumerative combinatorics}, 
        Vol. 2, Cambridge Studies in Advanced Mathematics, {\bf 62} 2000.  

\bibitem[To]{To} B. Totaro,
        {\em Euler and algebraic geometry},
        in print, Bulletin AMS 2007.

\bibitem[Tu]{Tu} H.L. Turrittin, 
        {\em The formal theory of systems of irregular homogeneous linear 
        difference and differential equations},
        Bol. Soc. Mat. Mexicana {\bf 5} (1960) 255--264. 

\bibitem[WP]{WP} J. Wimp and D. Zeilberger, 
        {\em Resurrecting the asymptotics of linear recurrences},
        J. Math. Anal. Appl.  {\bf 111}  (1985) 162--176.

\bibitem[WZ]{WZ} H. Wilf and D. Zeilberger,
        {\em An algorithmic proof theory for hypergeometric (ordinary and
        $q$) multisum/integral identities},
        Inventiones Math. {\bf 108} (1992)  575--633.

\bibitem[Ze]{Ze} D. Zeilberger,
        {\em A holonomic systems approach to special functions identities},
        J. Comput. Appl. Math. {\bf 32} (1990) 321--368.

\end{thebibliography}
\end{document}